\documentclass[11pt]{amsart}

\usepackage{amsmath,amssymb,amscd}
\usepackage[dvips]{graphics}

\newtheorem{theorem}{Theorem}[section]
\newtheorem{corollary}[theorem]{Corollary}

\newtheorem{lemma}[theorem]{Lemma}
\newtheorem{proposition}[theorem]{Proposition}
\newtheorem{remark}[theorem]{Remark}
\newtheorem{example}[theorem]{Example}

\def\PP{\mathbb{P}}
\def\ZZ{\mathbb{Z}}

\newcommand{\cE}{{\mathcal E}}

\newcommand{\cL}{{\mathcal L}}
\newcommand{\cO}{{\mathcal O}}
\newcommand{\cM}{{\mathcal M}}

\newcommand{\cF}{{\mathcal F}}

\newcommand{\cK}{{\mathcal K}}

\def\GL{\operatorname{GL}}
\def\pf{\operatorname{Pf}}

\def\coker{\operatorname{Coker}}
\def\M{\operatorname{Mat}}

\def\Id{\operatorname{Id}}
\def\pic{\operatorname{Pic}}
\def\rank{\operatorname{rank}}

\begin{document}

\title{Plane Curves as Pfaffians}
\author{Anita Buckley, Toma\v{z} Ko\v{s}ir}

\address{Department of Mathematics, University of Ljubljana,
Jadranska 19, 1000 Ljubljana, Slovenia  
{\rm e-mail:} {\it anita.buckley@fmf.uni-lj.si}}%

\begin{abstract}
Let $C$ be a smooth curve in $\PP^2$ given by an equation $F=0$ of degree $d$. 
In this paper we parametrise all linear pfaffian representations of $F$ by an open subset in the moduli space
$M_C(2,K_C)$.
We construct an explicit correspondence between pfaffian representations of $C$ and rank 2 vector bundles 
on $C$ with canonical determinant and no sections. 
\end{abstract}

\maketitle

\section{Introduction}
\label{introdef}

Let $k$ be an algebraically closed field and $C$ an irreducible curve defined by a polynomial $F$ 
of degree $d$ in $\PP^2$. 
When $C$ is smooth its genus equals 
$$g=\frac{1}{2}(d-1)(d-2).$$

These notes consider the following question. For a given homogeneous polynomial $F(x_0,x_1,x_2)$ of degree $d$ 
find a $2d\times 2d$ skew-symmetric matrix 
$$A= \left[\begin{array}{ccccc}
0       & L_{1\, 2} & L_{1\, 3}   & \cdots & L_{1\, 2d}\\
-L_{1\, 2} & 0        & L_{2\, 3} & \cdots & L_{2\, 2d}\\
-L_{1\, 3}  & -L_{2\, 3}& 0        &        &          \\
 \vdots   &  \vdots &          & \ddots &  \vdots         \\
 -L_{1\, 2d} & -L_{2\, 2d} &          &  \cdots  & 0
\end{array}\right]$$
with linear forms $L_{ij}=a_{ij}^0 x_0+a_{ij}^1 x_1+a_{ij}^2 x_2$
such that 
$$\pf A(x_0,x_1,x_2)=c F(x_0,x_1,x_2)\ \mbox{ for some }\ c\in k, c\neq 0.$$
Matrix $A$ is called \textit{linar pfaffian representation} of $C$. Its cokernel is a rank 2 vector bundle
on $C$. Throughout the paper we equate the notion of vector bundles and locally free sheaves.  

Two pfaffian representations $A$ and $A'$ are \textit{equivalent} if there exists $X\in\GL_{2d}(k)$ such that
$$A'=XAX^t.$$
A locally free sheaf $\mathcal{E}$ of rank 2 is \textit{stable} if for every invertible sheaf
$\cE\rightarrow \cF\rightarrow 0$ holds
$$\deg \cF > \frac{1}{2} \deg \cE.$$
Replacing $>$ by $\geq$ defines \textit{semistable}.

We find \textbf{all} linear pfaffian representations of $C$ (up to equivalence) and relate them to 
\textit{the moduli space} $M_C(2,K_C)$
of semistable rank $2$ vector bundles on $C$ with canonical determinant. 
An explicit construction of representations from the global sections of rank $2$
vector bundles with certain properties is given. 

In general the elements of $A$ can be homogenous polynomials of various degrees. Such pfaffian representations 
are considered in Theorem~\ref{thm1b}.
A good survey of the linear algebra of Pfaffians can be found in~\cite[Appendix D]{fulton}.
In~\cite[Chapter V. 2]{hartshorne} Hartshorne identifies the theory of rank 2 sheaves with the theory of 
ruled surfaces over $C$.

Study of pfaffian representations is strongly related to and motivated by determinantal representations. 
A \textit{linear determinantal representation} of $C$ is a $d\times d$ matrix of linear forms
$$M=x_0 M_0+x_1 M_1+x_2 M_2$$ 
satisfying
$$\det M=c F,$$
where $M_0, M_1, M_2\in\M_d(k)$ and $c\in k, c\neq 0$.
Two determinantal representations $M$ and $M'$ are \textit{equivalent} if there exist $X, Y\in\GL_{d}(k)$ such that
$$M'=XMY.$$
There are many more pfaffian than determinantal representations. Indeed, every determinantal representation $M$
induces \textit{ decomposable pfaffian representation}
$$\left[ \begin{array}{cc}
0 & M \\
-M^t & 0
\end{array}\right].$$ 
Note that the equivalence relation is well defined since
$$\left[ \begin{array}{cc}
0 & X M Y \\
-(X M Y)^t & 0
\end{array}\right]=\left[ \begin{array}{cc}
X & 0 \\
0 & Y^t
\end{array}\right]\left[ \begin{array}{cc}
0 & M \\
-M^t & 0
\end{array}\right]\left[ \begin{array}{cc}
X^t & 0 \\
0 & Y
\end{array}\right].$$ 

A brief outline of the paper is the following.
In Section~\ref{beginsec} we recall the parametrisation of linear determinantal representations by points on the 
Jacobian variety due to Vinnikov~\cite{vinnikov2}. We use similar ideas in 
Sections~\ref{secthree} and~\ref{secfour}
to parametrise all linear pfaffian representations by points in an open subset of 
the moduli space $M_C(2,K_C)$.
In the third section we give an explicit construction of the 
correspondence between pfaffian representations and vector bundles with certain properties.
In the fourth section we relate these vector bundles to the points on $M_C(2,K_C)$ not on the subvariety 
cut out by Cartier divisor $\Theta_{2,K_C}$. Section~\ref{decsec} considers decomposable bundles which
are parametrised by an open set in Kummer variety. Pfaffian representation arising from a decomposable bundle 
$\cF\oplus \left(\cF^{-1} \otimes \cO_C(d-1) \right)$ is exactly the decomposable representation
$$\left[ \begin{array}{cc}
0 & M \\
-M^t & 0
\end{array}\right],$$
where $M$ is the determinantal representation corresponding to $\cF\cong \coker M$.
Elliptic curves in Section~\ref{cubicsec} only allow decomposable representations.
In the last section we compute pfaffian representations of a genus 3 curve. 
To our knowledge there are not many results explicitly describing $M_C(2,K_C)$ for curves of higher genus.  
Methods considered in this paper could be helpful in 
finding all pfaffian representations and consequently the moduli space of $g>3$ curves. 

\section{Determinantal representations and the Jacobian} 
\label{beginsec}

There is a one to one correspondence between linear determinantal representations (up to equivalence) of $C$ and
line bundles (up to isomorphism) on $C$ with certain properties. This well known result is summed up in the 
following theorem of Beauville~\cite[Proposition 3.1]{beauville}.

\begin{theorem}
Let $C$ be a plane curve defined by a polynomial $F$ of degree $d$ and let $L$ be a line bundle of degree 
$\frac{1}{2}d(d-1)$ on $C$ with
$H^0(C, L(-1))=0$. Then there exists a $d\times d$ linear matrix $M$ with $\det M=F$ and an exact sequence 
\begin{eqnarray}\label{sxsqrk1}
0\rightarrow \bigoplus_{i=1}^{d} \cO_{\PP^2}(-1)\stackrel{M}{\longrightarrow} \bigoplus_{i=1}^{d} \cO_{\PP^2}
\rightarrow L \rightarrow 0.
\end{eqnarray}

Conversely, let $M$ be a linear $d\times d$ matrix with $\det M=F$. Then its cokernel is a line bundle of degree 
$\frac{1}{2}d(d-1)$ and
$H^0(C, \coker M(-1))=0$. 
\end{theorem}

Dolgachev explicitly described the above correspondence in~\cite{dolgachev}.Vinnikov also gave an explicit 
construction of the correspondence in ~\cite[Theorems 2-4]{vinnikov2}. Additionally he related 
determinantal representations to points on the Jacobian variety in the following way. 

It is known~\cite[Theorem 1.1]{cook} that representations $M$ and $M'$ are equivalent if and only if 
$\coker M$ and $\coker M'$ are isomorphic sheaves. Therefore the problem of classifying all linear representations 
of $F$ (up to equivalence) it the same to finding all line 
bundles with the property $\deg L = \frac{1}{2}d(d-1)$ and $H^0(C,L(-1))=0$. 
In order to simplify the notation we tensor the above by $\cO(1)$ and consider line bundles with
$\deg \cL = \frac{1}{2}d(d-1)-d = \frac{1}{2}d(d-3)$ and $H^0(C,\cL)=0$.
We call the bundles with this 
property \textit{nonexceptional} line bundles. Analogously, line 
bundles of degree $\frac{1}{2}d(d-3)$ and $H^0(C,\cL)\neq 0$ are called \textit{exceptional}.

\begin{lemma}
Let $C$ be of genus $g$. The exceptional divisor classes define a $g-1$ dimensional subvariety $W_{g-1}$ on 
the Jacobian variety $J$.
\end{lemma}

%\begin{proof}
%Recall that the Jacobian variety $J$ of $C$ is a $g=\frac{1}{2}(d-1)(d-2)$ dimensional Abelian variety. 
%For a fixed point $P_0\in C$ there is a regular map 
%$$\begin{array}{cccc} \mu: & C & \rightarrow & J\\
%                           & P_0 & \mapsto & 0,
%\end{array}$$
%which by additivity extends to all divisors on $C$. Moreover, for every $n$ it induces bijection
%$\mu:\, \mbox{Pic}^n(C)\rightarrow J$.
%
%Consider next the regular map
%$$\begin{array}{ccc}  C^{(n)} & \rightarrow & J\\
%                      \{ P_1,\ldots,P_n \} & \mapsto & \mu(P_1+\ldots +P_n),
%\end{array}$$
%where $C^{(n)}$ is the $n$th symmetric product of $C$. For $n\leq g$ the image of this map is an $n$ dimensional
%irreducible subvariety of $J$ which we denote $W_n$.
%
%Let $D$ be a divisor corresponding to $\cL$ with $H^0(X,\cL)\neq 0$. 
%Since 
%$$\deg \cL=\frac{1}{2}d(d-3)=g-1$$
%on $C$, there exist $P_1,\ldots,P_{g-1}$ on $C$ such that
%$$D\equiv P_1+\ldots +P_{g-1}.$$
%Clearly $\mu(D)=\mu(P_1+\ldots +P_{g-1})$. 
%
%Therefore all the exceptional line bundles define $W_{g-1}$, a g-1 dimensional subvariety on $J$.
%\end{proof}

This proves the following theorem. 

\begin{theorem}
All linear determinantal representations of $F$ 
(up to equivalence) can be parametrised by points on the Jacobian variety of $C$ not on the exceptional subvariety
$W_{g-1}$.
\end{theorem}

\section{Classification of pfaffian representations from the scratch}
\label{secthree}

In this section we consider the following problem. For a given homogeneous polynomial $F(x_0,x_1,x_2)$ of degree $d$ 
find all linear pfaffian representations. The main result of this section is an elementary proof of 
Theorem~\ref{propcondition} and explicit construction of representations from suitable vector bundles.\\

Following the ideas of Dolgachev~\cite{dolgachev} formulate the problem geometrically and coordinate free.
Let $U$ be a $2d$ dimensional vector space. It is well known that $\bigwedge^2 U$ can be identified 
with $2d\times 2d$ skew-symmetric matrices. Let $\Omega_k$ denote the set of vectors
$\sum_{i=1}^{k} v_i\wedge w_i$ in $\bigwedge^2 U$ where $\dim \{v_1,\ldots,v_k,w_1,\ldots,w_k \}=2k$. Elements of 
$\Omega_k$ are said to have \textit{irreducible length} $k$ since they can be written as a sum of $k$ and not 
less than $k$ pure nonzero products in $\bigwedge^2 U$.
In~\cite{westwick} it is shown that 
\begin{lemma}\label{westrk2}
$\Omega_k$ is isomorphic to the set of all rank $2k$ skew-symmetric matrices. 
\end{lemma}

\begin{proof}
Specifically, let 
$e_1,\ldots,e_{2d}$ be a basis for $U$ and $\{E_{ij}\}$ the standard basis for $2d\times 2d$ matrices.
Then the isomorphism equals 
$$ e_i\wedge e_j  \ \mapsto\   E_{ij}-E_{ji} $$
and extends linearly to
\begin{equation}
\label{explrk}
\sum_{i=1}^{2d} \alpha_i e_i\ \wedge\  \sum_{j=1}^{2d} \beta_j e_j \ \  \mapsto \ \ 
\sum_{i,j=1}^{2d}(\alpha_i \beta_j-\alpha_j \beta_i)(E_{ij}-E_{ji}). 
\end{equation}
\end{proof}

Next let $E$ be a 3 dimensional vector space and $\Phi$ a linear map
$$\Phi\!: \PP(E)\longrightarrow \PP(\bigwedge^2 U).$$ 
Note that $\Phi$ corresponds to a skew-symmetric matrix with linear forms as its elements. Alternatively, $\Phi$
is an element in $E^{\ast}\otimes (\bigwedge^2 U)$.

Let $\mathcal{P}_d\subset \PP(\bigwedge^2 U)$ be the hypersurface parametrising non-invertible linear maps. 
Choose a basis of $U$, then $\mathcal{P}_d$ is given by the pfaffian of a skew-symmetric matrix. The inverse image
of $\mathcal{P}_d$ under $\Phi$ is a plane curve of degree $d$ in $\PP(E)$.\\

Let $C$ be a smooth plane curve defined by $F$. Assume that $C$ admits a pfaffian representation $A$.

\begin{lemma}\label{rankpfaf}
For any $x\in C$ the corank of $A(x)$ equals 2.
\end{lemma}

\begin{proof}
Assume 
$$A=[a_{ij}^0 x_0+a_{ij}^1 x_1+a_{ij}^2 x_2]\ \mbox{ and }\ \pf A=c F,\ c\neq 0.$$
Denote by $\pf^{ij} A$ the pfaffian of the $(2d-2)\times(2d-2)$ skew-symmetric matrix obtained 
by removing the $i$th and $j$th rows and columns from $A$. Then
$$\frac{\partial F}{\partial x_k}(x)=\frac{1}{c} \sum_{i,j} a_{ij}^k \pf^{ij} A(x).$$
If for some $x\in C$ all $2d-2$ pfaffian minors vanish, then $x$ must be a singular point of $F$.
By our assumption $F$ is smooth, thus $\rank A(x) \geq 2d-2$ for all $x\in C$.
This ends the proof because rank of skew-symmetric matrices is even and $\det A=F^2=0$.
\end{proof}

Define the \textit{pfaffian adjoint} of $A$ to be the skew-symmetric matrix
\begin{equation*}\tilde{A}= \left[\begin{array}{ccl}
0 &       &         \\
  &\ddots &      (-1)^{i+j} \pf^{ij} A   \\
  &       &   \ddots                \\
  &       &  \ \ \ \ \ \ \ \ 0         
\end{array}\right].\end{equation*}
By analogy with determinants the following holds
\begin{equation} \label{pfaflinalg} 
\tilde{A} \cdot A=\pf A\cdot \mbox{Id}_{2d}.\end{equation}
More properties and linear algebra of Pfaffians can be found in~\cite[Appendix D]{fulton}.

By Lemma~\ref{rankpfaf} the cokernel of $A$ defines a $\rank 2$ bundle over $C$. Since $\pf A=c F$ the cokernel 
can be obtained from $\tilde{A}$ by (\ref{pfaflinalg}) . The properties of $\coker A$ are described in the following proposition.

\begin{proposition} \label{propofE} Let $A$ be a pfaffian representation of a smooth plane curve $C$ defined by a homogeneous 
polynomial $F$ of degree $d$. Then $\cE=\coker A$ is a rank 2 vector bundle on $C$ and
\begin{enumerate}
\item[(i)] $h^0(C,\cE)=2d$,
\item[(ii)] $H^0(C,\cE(-1))=H^1(C,\cE(-1))=0$,
\item[(iii)] $\det \cE = \bigwedge^2 \cE =\cO_C(d-1)$.
\end{enumerate}
\end{proposition}

\begin{proof}
By definition $\cE$ fits into the exact sequence
\begin{eqnarray}\label{sxsNo1}
0\rightarrow \bigoplus_{i=1}^{2d} \cO_{\PP^2}(-1)\stackrel{A}{\longrightarrow} \bigoplus_{i=1}^{2d} \cO_{\PP^2}
\rightarrow \cE \rightarrow 0.
\end{eqnarray}
Applying the functor $H^i(\PP^2,\ast)$ on (\ref{sxsNo1}) gives  
$$ \begin{array}{cccccccc}
H^0(\PP^2,\cO_{\PP^2}(-1)^{2d}) & \!\!\!\rightarrow \!\!\! & H^0(\PP^2,\cO_{\PP^2}^{2d}) & \!\!\!\rightarrow\!\!\! &
H^0(\PP^2,\cE) & \!\!\!\rightarrow\!\!\! & H^1(\PP^2,\cO_{\PP^2}(-1)^{2d}) & \!\!\!\rightarrow\cdots\!\!\! \\
  & &    & & & &   & \\
0 & & 2d & & & & 0 & 
\end{array}$$
where the bottom row denotes dimensions of the cohomology of projective space which is computed in
~\cite[Theorem III.5.1]{hartshorne}. Thus $\dim H^0(\PP^2,\cE)=2d $.

Next tensor (\ref{sxsNo1}) by $\cO_{\PP^2}(-1)$ and again apply the functor $H^i(\PP^2,\ast)$. 
This gives a long exact sequence with dimensions
$$ \begin{array}{cccccc}
H^0(\PP^2,\cO_{\PP^2}(-2)^{2d}) & \!\!\!\rightarrow \!\!\! & H^0(\PP^2,\cO_{\PP^2}(-1)^{2d}) 
& \!\!\!\rightarrow\!\!\! & H^0(\PP^2,\cE(-1)) & \!\!\!\rightarrow\!\!\! \\
\| & & \| & & & \\
0 & & 0 & & & \\
  & &   & & & \\
  H^1(\PP^2,\cO_{\PP^2}(-2)^{2d}) & \!\!\!\rightarrow \!\!\! & H^1(\PP^2,\cO_{\PP^2}(-1)^{2d}) 
  & \!\!\!\rightarrow\!\!\! & H^1(\PP^2,\cE(-1)) & \!\!\!\rightarrow\!\!\! \\
  \| & & \| & & & \\
  0 & & 0 & & & \\
  & &   & & & \\
  H^2(\PP^2,\cO_{\PP^2}(-2)^{2d}) & \!\!\!\rightarrow \!\!\! & \cdots 
    &  &  &  \\
    \| & &  & & & \\
  0 & &  & & & 
\end{array}$$
where $H^2(\PP^2,\cO_{\PP^2}(-2)^{2d})\cong H^0(\PP^2,\cO_{\PP^2}(-1)^{2d})=0$ by Serre duality. We proved that
$H^0(\PP^2,\cE(-1))=H^1(\PP^2,\cE(-1))=0$. 

Since $\cE$ is supported on $C$ we proved (i) and (ii).

In order to prove (iii) apply the functor $\mathcal{H}om_{\cO_{\PP^2}}(\ast,\cO_{\PP^2}(-1))$ to (\ref{sxsNo1}).
We get
$$ \begin{array}{cccccccc}
0 & \!\!\!\rightarrow \!\!\! & \mathcal{H}om(\cE,\cO_{\PP^2}(-1)) & \!\!\!\rightarrow \!\!\! 
& \mathcal{H}om(\cO_{\PP^2}^{2d},\cO_{\PP^2}(-1)) & \!\!\!\rightarrow \!\!\! 
& \mathcal{H}om(\cO_{\PP^2}(-1)^{2d},\cO_{\PP^2}(-1)) & \!\!\!\rightarrow \!\!\! \\
 & & & & & & & \\
  & \!\!\!\rightarrow \!\!\! & \mathcal{E}xt^1(\cE,\cO_{\PP^2}(-1)) & \!\!\!\rightarrow \!\!\! 
& \mathcal{E}xt^1(\cO_{\PP^2}^{2d},\cO_{\PP^2}(-1)) & \!\!\!\rightarrow \!\!\! 
&           \cdots                                    & ,
\end{array}$$
where 
\begin{eqnarray*}
\mathcal{H}om(\cO_{\PP^2}^{2d},\cO_{\PP^2}(-1))\cong (\cO_{\PP^2}^{2d})^{\vee}\otimes \cO_{\PP^2}(-1) 
\cong \cO_{\PP^2}(-1)^{2d} ,\\
\mathcal{H}om(\cO_{\PP^2}(-1)^{2d},\cO_{\PP^2}(-1))\cong (\cO_{\PP^2}(-1)^{2d})^{\vee}\otimes \cO_{\PP^2}(-1) 
\cong \cO_{\PP^2}^{2d} 
\end{eqnarray*}
and  $\mathcal{E}xt^1(\cO_{\PP^2}^{2d},\cO_{\PP^2}(-1))=0$ by~\cite[Proposition III.6.3]{hartshorne}.
This implies that $\mathcal{E}xt^i(\cE,\cO_{\PP^2}(-1))=0$ for $i\neq 1$ and the above sequence is isomorphic to
\begin{eqnarray}\label{sxsNo1trans}
0\rightarrow \bigoplus_{i=1}^{2d} \cO_{\PP^2}(-1)\stackrel{A^t}{\longrightarrow} \bigoplus_{i=1}^{2d} \cO_{\PP^2}
\rightarrow \mathcal{E}xt^1(\cE,\cO_{\PP^2}(-1))  \rightarrow 0.
\end{eqnarray}
Thus we obtain
\begin{eqnarray*}
\cE & \cong & \coker A \cong \coker A^t \cong \mathcal{E}xt^1(\cE,\cO_{\PP^2}(-1)) \cong \\
    &   & \mathcal{E}xt^1(\cE,\cO_{\PP^2}(-3)\otimes \cO_{\PP^2}(2)) \cong 
          \mathcal{E}xt^1(\cE,\cO_{\PP^2}(-3))\otimes \cO_{\PP^2}(2) \cong \\
    &   & \mathcal{H}om_C(\cE,\cO_{C}(d-3)) \otimes \cO_{\PP^2}(2) \cong \mathcal{H}om_C(\cE,\cO_{C}(d-1))
\cong \\
\cE^{\vee}\otimes \cO_{C}(d-1) & & 
\end{eqnarray*}
since by Serre duality 
$$\mathcal{E}xt^1(\cE,\cO_{\PP^2}(-3))\cong \mathcal{E}xt^1(\cE,\omega_{\PP^2})\cong
\mathcal{H}om_C (\cE, \omega_C)\cong \mathcal{H}om_C(\cE,\cO_{C}(d-3)).$$
Finally,
$$\wedge^2 \cE \cong \cE \wedge (\cE^{\vee}\otimes \cO_{C}(d-1))  \cong \cO_C(d-1).$$
\end{proof}

In the sequel the reverse problem will be considered. We will give an explicit construction of pfaffian 
representation from a vector bundle with certain properties.
\begin{theorem}\label{constrofrep}
Let $C$ be a smooth plane curve 
of degree $d$. To every rank 2 vector bundle $\cE$ on $C$ with properties 
\begin{enumerate}
\item[(i)] $h^0(C,\cE)=2d$,
\item[(ii)] $H^0(C,\cE(-1))=0$,
\item[(iii)] $\det \cE = \bigwedge^2 \cE =\cO_C(d-1)$
\end{enumerate}
we can assign a pfaffian representation $A_{\cE}$. 
In particular, isomorphic bundles induce equivalent representations.
\end{theorem}

\begin{proof}
Let $U=H^0(C,\cE)$ be the $2d$ dimensional vector space of global sections of $\cE$. We will define a map
$\psi$ from $C$ to the space of $2d\times 2d$ skew-symmetric matrices with entries from the space of homogeneous 
polynomials of degree $d-1$, such that $\psi^{-1}(\mathcal{P}_d)=C$.

Choose a basis $\{s_1,\ldots,s_{2d}\}$ for $U$ and define
$$C \ni x \ \stackrel{\psi}{\mapsto} \ \sum_{1\leq i<j\leq 2d}(s_i(x) \wedge s_j(x)) (E_{ij}-E_{ji})
= \left[\begin{array}{ccl}
0 &       &         \\
  &\ddots &      s_i(x)\wedge s_j(x)   \\
  &       &   \ddots                \\
  &       &  \ \ \ \ \ \ \ \ 0         
\end{array}\right].$$
Since $s_i\wedge s_j\in \bigwedge^2 U$, by property (iii) the map $\psi$ extends to
$$\Psi\!: \PP(E)\longrightarrow \PP(\bigwedge^2 U)$$ 
given by a linear system of plane curves of degree $d-1$. In other words, $\Psi$ is a tensor in 
$S^{d-1}(E^{\ast})\otimes (\bigwedge^2 U)$. In coordinates it equals to a
$2d\times 2d$ skew-symmetric matrix $B(x_0,x_1,x_2)$ with entries from the space of homogeneous 
polynomials of degree $d-1$. 

From the definition of $B$ and isomorphism (\ref{explrk}) in Lemma~\ref{westrk2} it follows that
at any point $x\in C\ \rank B(x)=2$.

Before we proceed, we prove that a different basis $\{s'_1,\ldots,s'_{2d}\}$ for $U$ induces equivalent representation
$B'$. 
Indeed, 
$$\left[\begin{array}{c}
s'_1\\
\vdots \\
s'_{2d}
\end{array}\right]=R \cdot
\left[\begin{array}{c}
s_1\\
\vdots \\
s_{2d}
\end{array}\right]$$
for some constant invertible matrix $R$, hence
\begin{equation*}B'= \left[\begin{array}{ccl}
0 &       &         \\
  &\ddots &      s'_i\wedge s'_j   \\
  &       &   \ddots                \\
  &       &  \ \ \ \ \ \ \ \ 0         
\end{array}\right]=R\cdot  \left[\begin{array}{ccl}
0 &       &         \\
  &\ddots &      s_i\wedge s_j   \\
  &       &   \ddots                \\
  &       &  \ \ \ \ \ \ \ \ 0         
\end{array}\right]\cdot R^t=R\cdot  B \cdot R^t.\end{equation*}

As before, let $F$ be the defining polynomial for $C$. Denote by $M$ a $4\times 4$ submatrix of $B$, obtained by
deleting $2d-4$ rows and columns with the same indexes. Then $F \,|\, \pf M$ since $\pf M(x)=0$ for all $x\in C$.
Consider next a $6\times 6$ skew-symmetric submatrix $N$ of $B$ and its pfaffian adjoint $\tilde{N}$. 
By (\ref{pfaflinalg}) we get
\begin{eqnarray*}
\pf N \cdot \pf \tilde{N}=(\pf N)^{\mbox{order} N/2},\\
\pf \tilde{N}=(\pf N)^{\mbox{order} N/2-1}=(\pf N)^2.
\end{eqnarray*}
The entries of $\tilde{N}$ are Pfaffians of $4\times 4$ submatrices, hence $F^3 \,|\, (\pf N)^2$.
Since $C$ is irreducible, $F^2 \,|\, \pf N$. By repeating this process we obtain that
$F^{d-2}$ divides all the Pfaffians of $(2d-2)\times (2d-2)$ skew-symmetric submatrices of $B$. These are exactly
$\pf^{ij} B$ defined in the proof of Lemma~\ref{rankpfaf}. This means that
$$A=\frac{1}{F^{d-2}}\, \tilde{B}$$
is a matrix with entries in $E^{\ast}$.
Since $\rank B(x)=2$ for all $x\in C$, we get $\rank A(x)=2d-2$. Therefore $\pf A$ is a hypersurface of degree
$d$ vanishing on $C$ unless $\pf A$ is identically zero. This implies $\pf A=c F$ for some constant $c$.

It remains to prove that $\pf A\neq 0$.
We will use the following remark to describe the sections of $\cE(-1)$.
\begin{remark}\label{sectminus} \rm{
Consider the map 
$$\begin{array}{cccc}
r: & C & \longrightarrow & \PP(U) \\
   & x & \mapsto         & [s_1(x),\ldots,s_{2d}(x)].
\end{array}$$
Every global section of $\cE$ can be seen as a hyperplane in $\PP(U)$. 
A nonzero element in $H^0(C,\cE(-1))$ is a section $s$ of $\cE$ such that the divisor
$$r(C)\ \cap\ \left( \mbox{the hyperplane in }\PP(U) \mbox{  defined by }s\right) $$
equals
$r(C\cap L)$ for some line $L$. All divisors on $C$ which are cut out by a 
line are equivalent. Hence for any line $L$ the divisor $r(C\cap L)$ on the image $r(C)$ is cut out by a hyperplane 
in $\PP(U)$.
}\end{remark}

Choose $L$ such that it intersects $C$ in $d$ distinct points $\{p_1,\ldots,p_d\}$.
Assumption (ii) and Remark~\ref{sectminus} imply that the divisor
$\{r(p_1),\ldots,r(p_d)\}$
%\begin{equation}\label{rankdivisor}
%\left\{
%\begin{array}{c}
%r(p_1)\\ 
%r(p_2)\\ 
%\vdots\\
%r(p_d)\end{array}
%\right\}
%=
%\left\{
%\begin{array}{cccc}
%s_1(p_1) & \cdots &  \cdots & s_{2d}(p_1) \\
%s_1(p_2) & \cdots &  \cdots &  s_{2d}(p_2) \\
%\vdots &  & & \vdots \\
%s_1(p_d) & \cdots &  \cdots &  s_{2d}(p_d) 
%\end{array}
%\right\}
%\end{equation}
on $r(C)$ is NOT 
cut out by a hyperplane in $\PP(U)$.
This proves that the images 
$$\{\Psi(p_1),\ldots,\Psi(p_d)\}$$
span a $d-1$ (projectively) 
dimensional subspace $W$ in $\PP(\bigwedge^2 U)$. Indeed, chose coordinates so that
the points $r(p_i)=[s_1(p_i), \ldots, s_{2d}(p_i)]$ are defined by 
$$\left[
\begin{array}{cccccccccc}
0 & \ldots & 1 & \ldots & 0 & 0 & \ldots & 0 & \ldots & 0 \\
0 & \ldots & 0 & \ldots & 0 & 0 & \ldots & 1 & \ldots & 0
\end{array}\right]$$
with 1 on the $i$th and $d+i$th position. Then $W$ corresponds to the space
of $2d\times 2d$ skew-symmetric matrices
$$\left[
\begin{array}{cc}
0 & \diamond \\
-\diamond & 0
\end{array}\right],$$ where $\diamond$ denotes diagonal $d\times d$ matrices.

By (iii) the image $\Psi(L)$ is a Veronese curve of degree $d-1$. Let $l_1,l_2$ be the local parameters of 
$L$, then $\Psi(L)=l_1^d+\ast\, l_1^{d-1}l_2+\cdots+\ast\, l_2^d$  is determined by $d$ points. 
This proves that $\Psi(L)$ lies in $W$. A general point $\Psi(p)\in\Psi(L)$ does not belong to any 
hyperplane in $W$ spanned by $d-1$ points $\Psi(p_i)$. Hence $\Psi(p)$ can be written as a linear combination of
the points $\Psi(p_i)$ with nonzero coefficients. In other words, $\Psi(p)$ can not be written shorter than
a sum of $d$ pure nonzero wedge products in $\PP(\bigwedge^2 U)$. By Lemma~\ref{westrk2}, 
$\Psi(p)\in \Omega_d$  corresponds to a skew-symmetric matrix of rank $2d$. Therefore 
$\pf B(p)\neq 0$ and hence $\pf A(p)\neq 0.$
\end{proof}

By the classic result of Cook~\cite[Theorem 1.1]{cook} two representations are equivalent if and only if 
the corresponding cokernels are isomorphic sheaves.
Together with the above considerations it implies

\begin{theorem}\label{propcondition}
There is a one to one correspondence between linear pfaffian representations of $F$ (up to equivalence) and
rank 2 bundles (up to isomorphism) on $C$ with the property 
\begin{eqnarray} \label{propertyE} \det \cE \cong \cO_C(d-1)\ \mbox{ and }\ H^0(C,\cE(-1))=0.\end{eqnarray}
\end{theorem}

\begin{proof}
It remains to show that (i) in Theorem~\ref{constrofrep} follows from conditions (ii) and (iii). 
For every rank 2 bundle $\cE\cong \cE^{\vee}\otimes (\bigwedge^2 \cE)$ holds. Applying (iii) gives
$\cE \cong \cE^{\vee}\otimes \cO_C(d-1)$. Then
\begin{eqnarray*}
H^1(C,\cE(-1)) \\
\cong H^1(C,\cE^{\vee}\otimes \cO_C(d-2)) \cong 
H^0(C,\cE\otimes \cO_C(2-d)\otimes \cO_C(d-3)) \\
\cong H^0(C,\cE(-1))
\end{eqnarray*}
by Serre duality. Hence (ii) implies $H^1(X,\cE(-1))=0$.

Let $L$ be a section of $\cO_C(1)$ and consider the exact sequence
\begin{equation}
0 \longrightarrow \cE(-1) \longrightarrow \cE\longrightarrow  \cE|_L \longrightarrow 0.
\end{equation}
Note that $\cE|_L=\cE\otimes \cO_L$ is supported on a set of points. Taking cohomology gives a long
exact sequence 
$$ \begin{array}{cccccc}
H^0(C,\cE(-1)) & \!\!\!\rightarrow \!\!\! & H^0(C,\cE) 
& \!\!\!\rightarrow\!\!\! & H^0(C,\cE|_L) & \!\!\!\rightarrow\!\!\! \\
\| & &  & &  & \\
0 & &  & &  & \\
  & &   & & & \\
  H^1(C,\cE(-1)) & \!\!\!\rightarrow \!\!\! & H^1(C,\cE) 
  & \!\!\!\rightarrow\!\!\! & H^1(C,\cE|_L) & \!\!\!\rightarrow\!\!\! \\
  \| & &  &  & \| & \\
  0 & &  &  & 0 & 
\end{array}$$
This proves that $H^1(C,\cE)=0$.

%Recall that~\cite[Corollary V.2.7]{hartshorne} any locally free sheaf of rank 2 on 
%a curve is an extension of invertible sheaves 
%\begin{equation*}
%0 \longrightarrow \cL \longrightarrow \cE\longrightarrow  \cM \longrightarrow 0.
%\end{equation*}
%By Riemann-Roch on $C$ we get
%\begin{eqnarray*}
%h^0(C,\cE)=h^0(C,\cE)-h^1(C,\cE)=\chi(\cE)=\chi(\cL)+\chi(\cM)\\
%=\deg \cL + 1-g+\deg \cM +1-g=\deg \cE+2-2g=d(d-1)+2-2g=2d.
%\end{eqnarray*}

It also proves $h^0(C,\cE)=2d$. Indeed, observe that 
$\cE|_L=\cE\otimes \cO_L$ is of rank 2 supported on the set $C\cap L$ of $d$ points, hence
$h^0(C,\cE|_L)=2d$.
\end{proof}

Theorem~\ref{propcondition} is a special case of the following Beauville's corollary, where representations
of hypersurfaces are studied via arithmetically Cohen-Macaulay (ACM) sheaves on $\PP^n$. One of the advantages of
our elementary proof is an explicit construction of representations from sheaves.

\begin{theorem}[Corollary 2.4 in \cite{beauville}] \label{thm1b}
Let $X$ be an integral hypersurface of degree $d$ in $\PP^n(x_0,\ldots,x_n)$ over a field $k$ with char $k\neq 2$.  
Moreover, let $E$ be an ACM vector bundle on $X$ of rank 2 with determinant $\cO_X(d+t)$. 
Then there exists a skew-symmetric matrix 
$A=(a_{ij})\in \mathbb{M}_l$ with $a_{ij}$ homogeneous polynomials in $x_0,\ldots,x_n$ of degree $d_i+d_j-t$ 
and an exact sequence
\begin{eqnarray}\label{sxsq}
0\rightarrow \bigoplus_{i=1}^{l} \cO_{\PP^n}(t-d_i)\stackrel{A}{\longrightarrow} \bigoplus_{i=1}^{l} \cO_{\PP^n}(d_i)
\rightarrow E \rightarrow 0,
\end{eqnarray}
where $X$ is defined by $\pf A=0$.

If $H^0(X,E(-1))=0$ and $t=-1$, the entries of $A$ are linear.
\end{theorem}

We conclude this section by Fujita's result~\cite[Example 6.4.16]{lazarsfeld}.

\begin{remark}\rm{Suppose that $C$ has genus $g\geq 2$. If $\cE$ is a bundle on $C$ such that $H^1(C,\cE)=0$, then
$\cE$ is ample. Being ample is equivalent to the condition that every quotient has strictly positive degree.
}\end{remark}

\section{The moduli space M(2,2(g-1))} \label{moduli}
\label{secfour}

In this section we relate the set of pfaffian representations with the moduli space
of semistable vector bundles. As described in Section~\ref{beginsec}, the idea arises from determinantal 
representations which can be parametrised by points on the Jacobian.

To any compact Riemann surface $X$ one can associate the pair $(JX,\Theta)$, the Jacobian and the Riemann theta function.
The geometry of the pair is strongly related to the geometry of $X$. This gives an idea that higher rank vector bundles 
define a non-abelian analogue of the Jacobian called \textit{moduli space} firstly due to the 
mathematicians of the Tata Institute. Much later physicists in Conformal Field Theory introduced pairs of 
moduli spaces and determinant line bundles on these moduli spaces.  
This has made a clear analogy with the Jacobian pair.\\

As before let $C$ be a smooth plane curve of degree $d$ and genus $g$. The existence and properties of 
$$M_C(r,n),$$ 
the moduli space of semistable vector bundles on $C$ of rank $r$ and degree $n$,
were established in~\cite{newstead}, ~\cite{seshadri} and more modern treatment can be found in~\cite{potier}. 
It is known that $M_C(r,n)$ is an irreducible, normal projective variety with an open subset 
$M_C^s(r,n)$ corresponding to stable bundles. 

If $C$ has genus $g\geq 2$ then $M_C^s(r,n)$ is not empty and 
its dimension is $r^2(g-1)+1$. The singular points of $M_C(r,n)$ are exactly $M_C(r,n)\backslash M_C^s(r,n)$.

If $C$ is an elliptic curve then $M_C^s(r,n)$ is empty.\\

One can restrict the study to the moduli space 
$$M_C(r,\cL)$$
of (semistable) rank $r$ vector bundles on $C$ with determinant $\cL$. As described 
in~\cite{beauvilleTh}, ~\cite{narasimhan} 
it is a closed subvariety in $M_C(r,\deg \cL)$. Moreover, $M_C^s(r,\cL)$ is a closed subvariety in $M_C^s(r,\deg \cL)$. 
The determinant can be fixed since the moduli space $M_C(r,\deg \cL)$ is, up to a finite \'{e}tale covering,
the product of $M_C(2,\cL)$ with the Jacobian $J C$. 

Drezed and Narasimhan~\cite{narasimhan} show that $\pic (M_C(r,\cL))\cong \ZZ$ is generated by geometrically 
defined Cartier divisors $\Theta_{r,\cL}$ in $M_C(r,\cL)$. For example, when $\deg \cL=r(g-1)$ then 
$$\chi(\cE)=0\ \mbox{ for all }\ \cE\in M_C(r,\cL)$$
and
$$\Theta_{r,\cL}=\left\{ \cE\in M_C(r,\cL)\ :\ h^0(C,\cE)\neq 0 \right\}$$
is naturally such a divisor.

\begin{theorem} \label{thmpfafmain} Let $C$ be a curve defined by a polynomial $F$ of degree $d$ in $\PP^2$.
There is a one to one correspondence between linear pfaffian representations of $F$ (up to equivalence) and
rank 2 bundles (up to isomorphism) on $C$ in the open set
$$M_C(2,\cO_C(d-3))\ \backslash \ \Theta_{2,\cO_C(d-3)}.$$
\end{theorem}

\begin{proof}
From now on we consider the moduli space 
$$M_C(2,\cO_C(d-3))$$
of (semistable) rank 2 vector bundles on $C$ with determinant $K_C\cong \cO_C(d-3)$, which is a closed subvariety in
$M_C(2,d(d-3))=M_C(2,2(g-1))$.

Recall from Theorem~\ref{propcondition} the one to one correspondence between the linear pfaffian representations 
of $C$ and rank 2 bundles on $C$ with property~(\ref{propertyE}) on page~\pageref{propertyE}. 
For the sake of clearer notation, after tensoring by $\cO_C(1)$, condition~(\ref{propertyE}) can be rewritten into 
\begin{eqnarray} \label{propertyEdpr} \det \cE \cong \cO_C(d-3)\ \mbox{ and }\ H^0(C,\cE)=0.\end{eqnarray}
Condition $h^0(C,\cE)=0$ implies that $\cE$ is semistable.
By the above considerations, bundles satisfying ~(\ref{propertyEdpr}) can be parametrised by the points on 
$$M_C(2,\cO_C(d-3))\ \backslash \ \Theta_{2,\cO_C(d-3)}.$$ 
This is an open subset of $M_C(2,\cO_C(d-3))$ since we cut out a Cartier divisor. 
\end{proof}

\section{Pfaffians arising from decomposable vector bundles}
\label{decsec}

In this section we find and explicitly describe linear pfaffian representations of $C$ (up to equivalence) 
arising from decomposable vector bundles. We parametrise them by an open set in the Kummer variety. 
Since decomposable vector bundles are never 
stable~\cite[Ex V.2.8]{hartshorne}, this open set lies in the singular locus of $M_C(2,\cO_C(d-3))$. The fact that
the moduli space of rank 2 bundles with canonical determinant  is singular along the Kummer variety can 
be found in~\cite{oxbury}.

\textit{Kummer variety} $\cK_C$ of $C$ is by definition the quotient of the Jacobian $J C$ by the involution 
$\cL \mapsto \cL^{-1}\otimes \cO_C(d-3)$.

\begin{figure}[ht] 
\begin{center}
\includegraphics*[1mm, 2cm][12cm, 9cm]{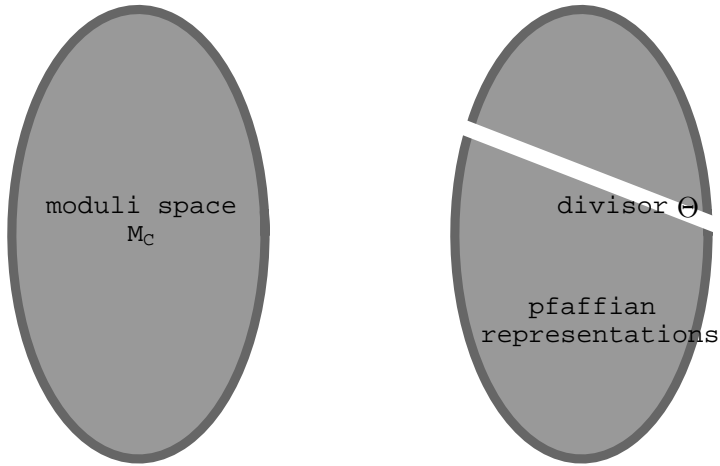}
  \put(-83,94){$\Theta$}
\caption{The darker areas represent the singular locus of $M_C$ and decomposable Pfaffians respectively}
\label{pict5}
\end{center}
\end{figure}

Recall that by Theorem~\ref{propcondition} and (\ref{propertyEdpr}) we need to find all decomposable rank 
2 bundles (up to isomorphism) with the property 
$$ \det \cE \cong \cO_C(d-3)\ \mbox{ and }\ H^0(C,\cE)=0.$$
Write $\cE=\cL \oplus \cM$. Then $\det \cE\cong \cL\otimes\cM \cong \cO_C(d-3)$ or equivalently 
$$\cM \cong \cL^{-1}\otimes \cO_C(d-3).$$ 
Observe that
$\cE \cong \cL \oplus \left(\cL^{-1}\otimes \cO_C(d-3)\right)$
has no sections if and only if $\cL$ and $\cL^{-1}\otimes \cO_C(d-3)$ have no sections. These can be 
calculated by the Riemann-Roch formula and Serre duality
$$h^0(C,\cL)-h^0\left(C,\cL^{-1}\otimes \cO_C(d-3)\right)=\deg \cL +1-g ,$$
since $\cO_C(d-3)$ is exactly the line bundle of the canonical divisor $K_C$.

This proves that every decomposable rank 2 bundles with the property (\ref{propertyEdpr}) is of the
form
$$\cE \cong \cL \oplus \left(\cL^{-1}\otimes \cO_C(d-3)\right),$$ 
where $\cL$ is a line bundle of degree $g-1=\frac{1}{2}d(d-3)$ with no sections.
These are exactly the \textit{nonexceptional} line bundles defined in
Section~\ref{beginsec}. We have seen that they correspond to the
points on the Jacobian variety of $C$ not on the exceptional subvariety $W_{g-1}$.

Conversely, by Section~\ref{beginsec} a point on the Jacobian variety induces a determinantal representation 
$M$ of $C$. Then 
$$\left[ \begin{array}{cc}
0 & M \\
-M^t & 0
\end{array}\right]$$
is a pfaffian representation of $C$ with decomposable cokernel.

Thus we proved
\begin{theorem}
There is a one to one correspondence between decomposable vector bundles in
$M_C(2,\cO_C(d-3))\ \backslash \ \Theta_{2,\cO_C(d-3)}$ and
the open subset of Kummer variety
$$\left(JC\ \backslash \ W_{g-1}\right)/ \equiv,$$
where $\equiv$ is the involution $\cL \mapsto \cL^{-1}\otimes \cO_C(d-3)$.
\end{theorem}

In the sequel explicitly construct the above correspondence from the sections of $\cE(1)$. As above,
$\cE = \cL \oplus \left(\cL^{-1}\otimes \cO_C(d-3)\right)$ with $\cL$ of degree $g-1$ and no sections.
Denote $\cL\otimes \cO_C(1)=\cL(1)$ by $\cF$. Then
\begin{equation*}
H^0(C,\cF(-1))=0\ \mbox{ and }\ \deg \cF=\frac{1}{2}d(d-1)=\deg \left(\cF^{-1}\otimes \cO_C(d-1)\right).
\end{equation*}
Let $H$ be a section of $\cO_C(1)$ and consider the exact sequence
\begin{equation*}
0 \longrightarrow \cF(-1) \longrightarrow \cF\longrightarrow  \cF|_H \longrightarrow 0,
\end{equation*}
where $\cF|_H$ is supported on a $d$ points $C\cap H$. Applying cohomology gives a long
exact sequence with dimensions
$$ \begin{array}{cccccccc}
H^0(C,\cF(-1)) & \!\!\!\rightarrow \!\!\! & H^0(C,\cF) 
& \!\!\!\rightarrow\!\!\! & H^0(C,\cF|_H) & \!\!\!\rightarrow\!\!\! & H^1(C,\cF(-1)) & \!\!\!\rightarrow, \!\!\! \\
0 & &  & & d & & 0 & 
\end{array}$$
since $h^0(C,\cF(-1))-h^1(C,\cF(-1))=0$ by the Riemann-Roch theorem. 
This proves that $h^0(C,\cF)=d$. The same way we show $h^0(C,\cF^{-1}\otimes \cO_C(d-1))=d$.
Let $\{f_1,\ldots,f_d\}$ and $\{m_1,\ldots,m_d\}$ be bases of the complete linear system of $\cF$ and 
$\cF^{-1}\otimes \cO_C(d-1)$ respectively. Then
$$\{(f_1,0),\ldots,(f_d,0),(0,m_1),\ldots,(0,m_d)\}$$ 
form a basis for $H^0(C,\cE(1))$. Obviously $(f_i,0)\wedge (0,m_j)=f_i\otimes m_j$,
whereas $(f_i,0)\wedge (f_j,0)$ and $(0,m_i)\wedge (0,m_j)$ equal 0. 
If we repeat the proof of Theorem~\ref{constrofrep}, the map $\Psi$
is of the form
$$\left[ \begin{array}{cc}
0 & f_i\otimes m_j \\
-f_j\otimes m_i & 0
\end{array}\right].$$
In coordinates it equals to a
$2d\times 2d$ skew-symmetric matrix $B$ with entries from the space of homogeneous 
polynomials of degree $d-1$ and zero diagonal blocks. As before
$$A=\frac{1}{F^{d-2}}\, \tilde{B}$$
is a pfaffian representation of $C$. In particular,
$$A=\left[ \begin{array}{cc}
0 & M \\
-M^t & 0
\end{array}\right],$$
where $M$ is the determinantal representation of $C$ corresponding to the nonexceptional line bundle $\cF$.

\section{Cubic curve} 
\label{cubicsec}

In this section $C$ will denote a curve defined by a smooth cubic polynomial $F$ in $\PP^2$.

\begin{corollary}\label{corcubint} On a cubic curve $C$ all linear pfaffian representations 
can be parametrised by the points on the Kummer variety $\cK_C-\{ \mbox{one point}  \}$. 
\end{corollary}

\begin{proof} Recall that on an elliptic curve $K_C\cong \cO_C$. Since $M_C^s(2,0)$ is empty, 
there are no stable bundles on $C$. On the other hand, by~\cite[\S 4]{beauvilleTh} the non-stable part of
$M_C(2,\cO_C)$ consists of decomposable vector bundles of the form $\cL\oplus\cL^{-1}$ for $\cL $ in the Jacobian $J C$.
But for $\cL\in J C$ the following conditions are equivalent:
\begin{itemize}
\item $h^0(C,\cL\oplus\cL^{-1})=0$,
\item $h^0(C,\cL)=0$,
\item $\cL \neq \cO_C$. 
\end{itemize}
Therefore
$$M_C(2,\cO_C)\ \backslash \ \Theta_{2,\cO_C}\ 
= \ \{ \cL\oplus\cL^{-1};\ \cL\in J C \}\ \backslash \ \{ \cO_C\oplus\cO_C \}.$$ 
\end{proof}

Vinnikov~\cite{vinnikov2} found an explicit one to one correspondence between the linear determinantal representations 
(up to equivalence) of $C$ and the points on an affine piece of $C$:

\begin{lemma}[\cite{vinnikov2}] \label{vinncubb} Every smooth cubic can be brought into the Weierstrass form
$$F(x_0,x_1,x_2)=-x_1 x_2^2+x_0^3+\alpha x_0 x_1^2+\beta x_1^3.$$
A complete set of determinantal representations of $F$ is
$$x_0 \Id+x_2 \left(\begin{array}{ccc} 0&1&0\\ 0&0&1\\ 0&0&0 \end{array}\right) 
+x_1 \left(\begin{array}{ccc} \frac{s}{2}&l&\alpha+\frac{3}{4}s^2\\ 0&-s&-l\\ -1&0&\frac{s}{2} \end{array}\right),$$
where $l^2=s^3+\alpha s+\beta.$ Note that the last equation is exactly the affine part $F(s,1,l)$.
\end{lemma}

Recall that the Jacobian of a cubic curve $C$ with $g=1$ is the curve itself and $J-\{W_0\}$ is an affine piece 
of $C$. In particular, Corollary~\ref{corcubint} implies that the complete set of pfaffian representations of $F$
(put in the Weierstrass form) equals
$$\left[ \begin{array}{cc}
0 & M \\
-M^t & 0
\end{array}\right],$$
where $M$ are the determinantal representations in Lemma~\ref{vinncubb}.

%We finish by listing all pfaffian representations of a given $F$.  
%By Theorem~\ref{thm1b} there are three possibilities for the sizes of of a skew-symmetric $M=(m_{ij})$:
%\begin{itemize}
%\item $M\in \mathbb{M}_6$ with all $m_{ij}$ linear, or the entries of $M$ have different degrees $\{0,1,2,3\}$, 
%for example 
%$$\left[ \begin{array}{cccccc}  &2&2&2&2&1 \\  & &1&1&1&0\\ & & & 1&1&0\\ & & & &1&0\\  & & & & &0\\
% & & & & &\end{array} \right]\ \mbox{ or }\ 
% \left[ \begin{array}{cccccc}  &2&2&3&2&2 \\  & &0&1&0&0\\ & & & 1&0&0\\ & & & &1&1\\  & & & & &0\\
% & & & & &\end{array} \right]$$
%\item $M\in \mathbb{M}_4$ and its entries have degrees 
%$$\left[ \begin{array}{cccc}  &2&2&2 \\  & &1&1\\ & & & 1\\ & & & \end{array} \right],\ 
%\left[ \begin{array}{cccc}  &1&1&1 \\  & &2&2\\ & & & 2\\ & & & \end{array} \right],\ 
%\left[ \begin{array}{cccc}  &3&3&3 \\  & &0&0\\ & & & 0\\ & & & \end{array} \right]\ \mbox{ or }\
%\left[ \begin{array}{cccc}  &0&0&0 \\  & &3&3\\ & & & 3\\ & & & \end{array} \right],$$
%\item $M=\left[ \begin{array}{cc}  0&F \\  -F&0 \end{array} \right].$
%\end{itemize}

%The corresponding cokernels are ACM and torsion free.
%Note that $M$ in Theorem~\ref{thm1b} is mimimal (i.e., if $\deg m_{ij}=0$ then $m_{ij}=0$) and thus unique.

\section{Examples of higher genus}
\label{seclast}

We start by an example of a genus 3 curve.

\begin{example}\rm{
Any non hyperelliptic curve $C$ of genus 3 is isomorphic to a plane quartic. In this case 
$M_C(2,\cO_C(1))\cong M_C(2,\cO_C)$
embeds as a Coble quartic hypersurface in $\PP^7$ and is singular along the Kummer variety $\cK_C$. For references 
check~\cite{narasimhanram},~\cite{laszlo},~\cite{beauvilleTh}.

For a given plane quartic Vanhaecke~\cite{vanhaeck} gives explicit equations of the Coble quartic hypersurface.  
First he finds the equations of the Kummer variety which represents the singular locus of the moduli space, 
from here the Coble quartic is obtained by integration. 
For example, the moduli space of 
$$C:\ x^4 - y z^3 - y^4=0$$ 
is the Coble hypersurface in $\PP^7$ with equation
\begin{eqnarray*}z_0^4 - z_2^4 - 2 z_1 z_2^2 z_3 - 
    z_1^2 z_3^2 - z_0 z_3^3 - 4 z_0 z_2^2 z_4 - 2 z_0^2 z_4^2 - 
    3 z_0 z_2 z_4^2 + z_4^4 - 4 z_0 z_1 z_2 z_5 \\
    + z_3 z_4 z_5^2 
    + z_2 z_5^3 - 4 z_0 z_1^2 z_6 
    - 3 z_0^2 z_3 z_6 - z_3 z_4^2 z_6 
    - 2 z_2 z_4 z_5 z_6 + z_1 z_5^2 z_6 - z_0 z_6^3 - 2 z_0^2 z_1 z_7 \\ - z_3^2 z_4 z_7 + 
    2 z_1 z_4^2 z_7 
    - z_2 z_3 z_5 z_7 + z_0 z_5^2 z_7 
    - z_2^2 z_6 z_7 - z_1 z_3 z_6 z_7 - 4 z_0 z_4 z_6 z_7 + z_1^2 z_7^2.
    \end{eqnarray*}
    
In the sequel we outline an algorithm for finding all pfaffian representations of the given $C$ (up to equivalence)
based on canonical forms of matrix pairs. 
This is a generalisation of Vinnikov's construction of determinant representations~\cite{vinnikov2}.
By Lancaster and Rodman~\cite[Theorem 5.1]{rodman} every pfaffian representation of $C$ can be put into
the skew-symmetric canonical form
$$
\substack{x\left(\!\!\!\!\begin{array}{cccccccc} 
\substack{0} \!\!&\!\! \substack{0} \!\!&\!\!  \substack{0} \!\!&\!\!  \substack{0} \!\!&\!\!  
\substack{0} \!\!&\!\! \substack{0} \!\!&\!\!  \substack{0} \!\!&\!\!  \substack{1} \\
             \!\!&\!\! \substack{0} \!\!&\!\!  \substack{0} \!\!&\!\!  \substack{0} \!\!&\!\! 
\substack{0} \!\!&\!\! \substack{0} \!\!&\!\!  \substack{1} \!&\!\!\!  \substack{0} \\
             \!\!&\!\!              \!\!&\!\!  \substack{0} \!\!&\!\!  \substack{0} \!\!&\!\!  
\substack{0} \!\!&\!\! \substack{1} \ \!&\!\!  \substack{0} \!\!&\!\!  \substack{0} \\
             \!\!&\!\!              \!\!&\!\!               \!\!&\!\!  \substack{0} \!\!&\!\!  
\substack{1} \!\!&\!\! \substack{0} \!\!&\!\!  \substack{0} \!\!&\!\!  \substack{0} \\
             \!\!&\!\!              \!\!&\!\!               \!\!&\!\!               \!\!&\!\! 
\substack{0} \!\!&\!\! \substack{0} \!\!&\!\!  \substack{0} \!\!&\!\!  \substack{0} \\
             \!\!&\!\!              \!\!&\!\!               \!\!&\!\!               \!\!&\!\!  
             \!\!&\!\! \substack{0} \!\!&\!\!  \substack{0} \!\!&\!\!  \substack{0} \\
             \!\!&\!\!              \!\!&\!\!               \!\!&\!\!               \!\!&\!\! 
             \!\!&\!\!              \!\!&\!\!  \substack{0} \!\!&\!\!  \substack{0} \\
             \!\!&\!\!              \!\!&\!\!               \!\!&\!\!               \!\!&\!\!  
             \!\!&\!\!              \!\!&\!\!               \!\!&\!\!  \substack{0} 
\end{array}\!\!\!\! \right) }+
\substack{z\left(\!\!\!\!\begin{array}{cccccccc} 
\substack{0} \!\!&\!\! \substack{0} \!\!&\!\!  \substack{0} \!\!&\!\!  \substack{0} \!\!&\!\!  
\substack{0} \!\!&\!\! \substack{0} \!\!&\!\!  \substack{1} \!\!&\!\!  \substack{0} \\
             \!\!&\!\! \substack{0} \!\!&\!\!  \substack{0} \!\!&\!\!  \substack{0} \!\!&\!\! 
\substack{0} \!\!&\!\! \substack{1} \!\!&\!\!  \substack{0} \!&\!\!\!  \substack{0} \\
             \!\!&\!\!              \!\!&\!\!  \substack{0} \!\!&\!\!  \substack{0} \!\!&\!\!  
\substack{1} \!\!&\!\! \substack{0} \ \!&\!\!  \substack{0} \!\!&\!\!  \substack{0} \\
             \!\!&\!\!              \!\!&\!\!               \!\!&\!\!  \substack{0} \!\!&\!\!  
\substack{0} \!\!&\!\! \substack{0} \!\!&\!\!  \substack{0} \!\!&\!\!  \substack{0} \\
             \!\!&\!\!              \!\!&\!\!               \!\!&\!\!               \!\!&\!\! 
\substack{0} \!\!&\!\! \substack{0} \!\!&\!\!  \substack{0} \!\!&\!\!  \substack{0} \\
             \!\!&\!\!              \!\!&\!\!               \!\!&\!\!               \!\!&\!\!  
             \!\!&\!\! \substack{0} \!\!&\!\!  \substack{0} \!\!&\!\!  \substack{0} \\
             \!\!&\!\!              \!\!&\!\!               \!\!&\!\!               \!\!&\!\! 
             \!\!&\!\!              \!\!&\!\!  \substack{0} \!\!&\!\!  \substack{0} \\
             \!\!&\!\!              \!\!&\!\!               \!\!&\!\!               \!\!&\!\!  
             \!\!&\!\!              \!\!&\!\!               \!\!&\!\!  \substack{0} 
\end{array}\!\!\!\! \right) }+
\substack{y\substack{\left( \!\!\!\! \begin{array}{cccccccc} 
0 \!\!&\!\! c_{12} \!\!&\!\! c_{13} \!\!&\!\! c_{14} \!\!&\!\! c_{15} \!\!&\!\! c_{16} \!\!&\!\! 
c_{17} \!\!&\!\! c_{18} \\
\!\!&\!\! 0 \!\!&\!\! c_{23} \!\!&\!\! c_{24} \!\!&\!\! c_{25} \!\!&\!\! c_{26} \!\!&\!\! 
c_{27} \!\!&\!\! c_{28} \\
\!\!&\!\!  \!\!&\!\! 0 \!\!&\!\! c_{34} \!\!&\!\! c_{35} \!\!&\!\! c_{36} \!\!&\!\! 
c_{37} \!\!&\!\! c_{38} \\
 \!\!&\!\!  \!\!&\!\!  \!\!&\!\! 0 \!\!&\!\! c_{45} \!\!&\!\! c_{46} \!\!&\!\! 
c_{47} \!\!&\!\! c_{48} \\
 \!\!&\!\!  \!\!&\!\!  \!\!&\!\!  \!\!&\!\! 0 \!\!&\!\! c_{56} \!\!&\!\! 
c_{57} \!\!&\!\! c_{58} \\
 \!\!&\!\!  \!\!&\!\!  \!\!&\!\!  \!\!&\!\!  \!\!&\!\! 0 \!\!&\!\! 
c_{67} \!\!&\!\! c_{68} \\
 \!\!&\!\!  \!\!&\!\!  \!\!&\!\!  \!\!&\!\!  \!\!&\!\!  \!\!&\!\! 
0 \!\!&\!\! c_{78} \\
 \!\!&\!\!  \!\!&\!\!  \!\!&\!\!  \!\!&\!\!  \!\!&\!\!  \!\!&\!\!  \!\!&\!\! 0 
\end{array}\!\!\!\! \right), }}
$$
which we denote by $A$. Since $\pf A$ equals the equation of $C$ we get
\begin{eqnarray*}\substack{
c_{48} } &\!\!\!\!\substack{=}\!\!\!\! & \substack{1},\\
\substack{ c_{15} }&\!\!\!\!\substack{=}\!\!\!\!& 
\substack{-c_{18}^2 + 2 c_{26} c_{28} + c_{28} c_{35} - c_{18} c_{36} - c_{36}^2 + 2 c_{25} 
      c_{38} + 2 c_{18} c_{28} c_{38} - c_{28} c_{36} c_{38} - c_{26} c_{38}^2 - 2 c_{35} 
        c_{38}^2 - c_{18} c_{38}^3}\\
  &\!\!\!\!\!\!\!\!&  \substack{ - c_{18} c_{45} - c_{36} c_{45} - c_{28} c_{38} c_{45} + 
      c_{38}^3 c_{45} - c_{45}^2 + c_{26} c_{46} - c_{28}^2 c_{46} - c_{35} c_{46} + c_{18} 
      c_{38} c_{46} + c_{36} c_{38} c_{46}+ c_{28} c_{38}^2 c_{46} }\\
   &\!\!\!\!\!\!\!\!&   \substack{ + 2 c_{38} c_{45} 
      c_{46} + c_{34} c_{56} + c_{24} c_{57} - c_{34} c_{38} c_{57} +
       c_{14} c_{58} - c_{28} c_{34} c_{58} - c_{24} c_{38} c_{58} + c_{34} c_{38}^2 c_{58} + 
      c_{23} c_{67}- c_{24} c_{38} c_{67}  }\\
    &\!\!\!\!\!\!\!\!&  \substack{+ c_{13} c_{68} - c_{24} c_{28} c_{68} - c_{14} c_{38} c_{68} -
       c_{23} c_{38} c_{68} + c_{24} 
      c_{38}^2 c_{68} + c_{12} c_{78} - c_{14} c_{28} c_{78} - c_{13} c_{38} c_{78} + c_{14} c_{38}^2 c_{78}, }\\
\substack{c_{16} }&\!\!\!\!\substack{=}\!\!\!\!& 
\substack{-c_{25} - c_{18} c_{28} + c_{28} c_{36} + c_{26} c_{38} + 2 c_{35} c_{38} + 
      c_{18} c_{38}^2 - c_{38}^2 c_{45} - c_{18} c_{46} - c_{36} 
      c_{46}- c_{28} c_{38} c_{46}} \\
    &\!\!\!\!\!\!\!\!&  \substack{ - 2 c_{45} c_{46} + c_{34} c_{57} + c_{24} c_{58} - c_{34} c_{38} c_{58} 
+ c_{24} c_{67} + 
      c_{14} c_{68} + c_{23} c_{68} - c_{24} c_{38} c_{68} + c_{13} c_{78} - c_{14} c_{38} c_{78}, }\\
\substack{c_{17}} &\!\!\!\!\substack{=}\!\!\!\!& 
\substack{-c_{26} - c_{35} - c_{18} c_{38} + c_{38} c_{45} + c_{28} c_{46} + c_{34} c_{58} + c_{24} c_{68} + c_{14} 
        c_{78},}\\
\substack{c_{37}} &\!\!\!\!\substack{=}\!\!\!\!& \substack{-c_{28} - c_{46},}\\
\substack{c_{27}} &\!\!\!\!\substack{=}\!\!\!\!& \substack{-c_{18} - c_{36} - c_{45},}\\
\substack{c_{47}} &\!\!\!\!\substack{=}\!\!\!\!& \substack{ -c_{38}.}
       \end{eqnarray*} 
There are 21 parameters $c_{ij}$ left in the representation. We eliminate them by acting on $A$ with an invertible 
constant matrix $P$ via
$$P\cdot A\cdot P^t.$$
It is easy to check that $P$ preserves the canonical form of $A$ if and only if 
$$P=\left[\begin{array}{cc}
Y^{-1} & Y^{-1} S \\
Y (T-S)^{-1} & Y (\Id+(T-S)^{-1}S)
\end{array}\right],$$
where
$$\substack{Y=\left(\begin{array}{cccc}
p&q&s&t \\
0&p&q&s\\
0&0&p&q\\
0&0&0&p
\end{array}\right),\ 
S=\left(\begin{array}{cccc}
a&b&c&d \\
0&a&b&c\\
0&0&a&b\\
0&0&0&a
\end{array}\right),\
T=\left(\begin{array}{cccc}
e&f&g&h \\
0&e&f&g\\
0&0&e&f\\
0&0&0&e
\end{array}\right).}$$
This enables us to reduce to 12 parameters $c_{ij}$
$$\substack{\left[\begin{array}{cccccccc}
0 & c_{12} & c_{13} & 0 & 
\substack{-c_{36}^2 - c_{36} c_{45} - c_{45}^2 + c_{26} c_{46} - c_{35} c_{46} + c_{23} c_{67}} & 
\substack{-c_{25} - c_{36} c_{46} - 2 c_{45} c_{46}} & 
\substack{-c_{26} - c_{35} } & 0 \\
  & 0 & c_{23} & 0 & c_{25} & c_{26} & \substack{-c_{36} - c_{45} } &  0\\
  &   & 0   & 0 & c_{35} & c_{36} & \substack{-c_{46}}       & 0        \\
  &   &     & 0 & c_{45} & c_{46} & 0          & 1\\
  &   &     &   & 0 & c_{56} & c_{57} & 0\\ 
  &   &     &   &   & 0   & c_{67} & 0\\
  &   &     &   &   &     & 0 & 0\\
  &   &     &   &   &     &   & 0
\end{array}\right]}$$    

Another action by $P$ which preserves the 0 and 1 elements in the above matrix yields
\begin{eqnarray*}
c_{12} & \mapsto & \substack{\frac{c_{12}+a(-2 c_{25} - 2 c_{36} c_{46} - 3 c_{45} c_{46} + a c_{56} + a c_{46} c_{67}) }{p^2} }\\
c_{13} & \mapsto & \substack{\frac{c_{13}+a(-c_{26} - 2 c_{35} - c_{46}^2 + a c_{57})}{p^2} } \\
c_{23} & \mapsto & \substack{\frac{c_{23}+a(-2 c_{36} - c_{45} + a c_{67})}{p^2} }\\
c_{25} & \mapsto & \substack{\frac{c_{12} - a c_{25} - c_{23} c_{46} - 2 a c_{45} c_{46} - e c_{25} + a e c_{56}}{a-e} }\\
c_{26} & \mapsto & c_{26}\\
c_{35} & \mapsto & \substack{\frac{c_{13} - e c_{35} - a(c_{26} + c_{35} + c_{46}^2 - e c_{57})}{a-e}  }\\
c_{36} & \mapsto & \substack{\frac{c_{23} - e c_{36} - a(c_{36} + c_{45} - e c_{67})}{a-e} }\\
c_{45} & \mapsto & c_{45}\\
c_{46} & \mapsto & c_{46}\\
c_{56} & \mapsto & \substack{\frac{(c_{12} - c_{23} c_{46} + e(-2 c_{25} - 2 c_{45} c_{46} + e c_{56})) p^2}{(a-e)^2} }\\
c_{57} & \mapsto & \substack{\frac{(c_{13} + e(-c_{26} - 2 c_{35} - c_{46}^2 + e c_{57})) p^2}{(a-e)^2} }\\
c_{67} & \mapsto & \substack{\frac{(c_{23} + e(-2 c_{36} - c_{45} + e c_{67})) p^2}{(a-e)^2} }\\
 \end{eqnarray*}
Observe that 
$c_{13} c_{57} - c_{35}(c_{35} + c_{26} + c_{46}^2), c_{23} c_{67} - c_{36}(c_{36} + c_{45})$
are invariant under this action.

We are left with 12 parameters $c_{ij}$ modulo $P(a, e, p)$ action. The condition
$\pf A=x^4 - y z^3 - y^4$
gives three additional equations
\begin{eqnarray*}
\substack{2 c_{36}^2 c_{45} + 2 c_{36} c_{45}^2 +  c_{45}^3 + 2 c_{35} c_{45} c_{46} +
c_{25} (c_{26} + 2 c_{35} + c_{46}^2) + c_{23} c_{46} c_{57} } & \substack{= } & \\
 \substack{     2 c_{26} c_{45} c_{46} + c_{13} c_{56} + c_{12} c_{57} + 2 c_{23} c_{45} c_{67}, }  & & \\ 
\substack{  } & &\\
\substack{c_{36}^4 + 2 c_{36}^3 c_{45} + c_{25}^2 c_{46} + 2 c_{25} c_{45} c_{46}^2 + c_{26}^2 (-c_{35} 
+ c_{46}^2) + c_{13} c_{45} c_{56} + c_{23} c_{25} c_{57}  + 2 c_{23} c_{45} c_{46} c_{57} + c_{12} c_{35} c_{67} + } & &\\
\substack{  c_{23}^2 c_{67}^2 
+ 2 c_{36}^2 (c_{45}^2 - c_{23} c_{67}) + c_{36}(c_{45}^3 + c_{25} c_{46}^2 + c_{13} c_{56} - c_{12} c_{57} 
+ c_{23} c_{46} c_{57} - 2 c_{45} (c_{35} c_{46} + c_{23} c_{67})) } & \substack{= } & \\
\substack{1 + c_{25} c_{35} c_{45} + 2 c_{35} c_{45}^2 c_{46} + c_{23} c_{35} c_{56} + c_{12} c_{46} c_{56} 
+ c_{13} c_{25} c_{67}+ c_{23} c_{45}^2 c_{67} + c_{23} c_{35} c_{46} c_{67} + } & &\\
\substack{ c_{26} (c_{35}^2 - c_{25} c_{36} + 2 c_{36}^2 c_{46} + 2 c_{36} c_{45} c_{46} + c_{45}^2 c_{46} 
+ c_{35} c_{46}^2 + c_{23} c_{56} - c_{13} c_{57} - 2 c_{23} c_{46} c_{67}), }& & \\
\substack{  } & &\\
\substack{      c_{26}^2 + c_{26} c_{35} + c_{35}^2 + 2 c_{25} c_{36} + c_{25} c_{45} + 2 c_{36}^2 c_{46} + 4 c_{36} c_{45} c_{46} + 3 c_{45}^2 c_{46} + c_{35} c_{46}^2 }& \substack{= } &\\
\substack{            c_{26} c_{46}^2 + c_{13} c_{57} + c_{12} c_{67} + c_{23} (c_{56} + c_{46} c_{67}).} & &
            \end{eqnarray*}
          
We proved that all pfaffian representations of the given quartic curve form a 

\noindent $6=12 - 3$(parameters $a,e,p$) $- 3$(relations among $c_{ij}$) 
dimensional affine variety. 
}\end{example}

%\begin{remark}\rm{
Let $C$ be a generic plane curve of genus $g\geq 3$. By~\cite{oxbury} the moduli space $M_C(2,\cO_C(d-3))$ 
of rank 2 bundles with canonical determinant embeds into $|2 \Theta |$. 
Its singular locus is isomorphic to the Kummer variety. 
Moreover, the embedding restricts to the Kummer map on the singular locus. 

There are not many results in the literature explicitly describing the above moduli spaces for curves of genus $g>3$. 
Finding all pfaffian representations of a given curve $C$ would provide a description of the open set
$M_C(2,\cO_C(d-3))\ \backslash \ \Theta_{2,\cO_C(d-3)}.$ 

Vanhaecke's method could be used for curves of higher genus. 
Start by finding all determinantal representations of $C$ (which are smaller than Pfaffians)
using Vinnikov's method of canonical forms~\cite{vinni1}.
Determinantal representations induce pfaffian representations of the form 
$$\left[ \begin{array}{cc}
0 & \diamond \\
- \diamond^t & 0
\end{array}\right],$$
which by Section~\ref{decsec} correspond to decomposable bundles. These define the singular locus of the moduli space.
The equations of the moduli space could be then found by integration.

%}\end{remark}

\end{document}